\documentclass[11pt]{article}

\usepackage{graphicx}
\usepackage{amsmath}
\usepackage{amsthm}
\usepackage{amssymb}
\usepackage[mathscr]{eucal}

\title{The Grothendieck-Katz Conjecture for certain locally symmetric varieties}
\author{Benson Farb and Mark Kisin \thanks{Both authors are
supported in part by the NSF.}}

\theoremstyle{plain}
\newtheorem{theorem}{Theorem}

\newtheorem{conjecture}[theorem]{Conjecture}

\def\proof{{\bf {\medskip}{\noindent}Proof. }}

\def\title{\em}

\def\bar{\overline}

\newcommand\G{\mbox{\bf G}}
\newcommand\R{{\bf R}}

\newcommand\C{{\bf C}}

\newcommand\Z{\mbox{\bf Z}}

\newcommand\Q{{\bf Q}}
\newcommand\gal{{\rm gal}}
\newcommand\tor{{\rm tor}}
\newcommand\alg{{\rm alg}}
\newcommand\End{\mbox{End}}
\newcommand\Der{\mbox{Der}}

\DeclareMathOperator\ad{ad}
\DeclareMathOperator\Lie{Lie}

\DeclareMathOperator\der{der}

\DeclareMathOperator\SL{SL}

\DeclareMathOperator\GL{GL}

\DeclareMathOperator\Sp{Sp}

\begin{document}
\maketitle

\begin{abstract}
Using Margulis's results on lattices in semisimple Lie groups, we prove the Grothendieck-Katz $p$-Curvature Conjecture 
for certain locally symmetric varieties, 
including the moduli space of abelian varieties ${\cal A}_g$ when $g > 1.$
\end{abstract}


\section{Introduction}
In this paper we prove certain cases of the 
well-known $p$-curvature conjecture of Grothendieck-Katz stated below.  The conjecture posits that one can deduce algebraic solutions of certain differential equations when one has solutions after reducing modulo a prime for almost every prime.   Our main purpose is to point out that rigidity theorems from the theory of discrete subgroups of Lie groups can be fruitfully applied to this problem.

More precisely, let $X$ be a smooth connected variety over $\C$ and $V$ a vector bundle on $X$ equipped with an integrable connection 
$\nabla.$ Then there is a finitely generated $\Z$-algebra $R \subset \C$ such that $X$ arises from a smooth $R$-scheme and 
$(V,\nabla)$ descends to a vector bundle with integrable connection on this $R$-scheme. We will again denote by $X$ and $(V,\nabla)$ 
the corresponding objects over $R.$

For any maximal ideal $\mathfrak p$ of $R$, we can reduce ${\rm mod\,} {\mathfrak p}$ to obtain 
a differential equation $(V/{\mathfrak p}V,\nabla)$ on a smooth scheme over a finite field of characteristic $p>0$.  Attached to this system is an invariant, the 
$p$-curvature of $(V/{\mathfrak p}V,\nabla)$, which is an ${\cal O}_X$-linear map 
$$ \psi_p(V,\nabla): \underline{\Der}({\cal O}_X\otimes R/\mathfrak p, {\cal O}_X\otimes R/\mathfrak p) \rightarrow \underline\End_{{\cal O}_X}(V/\mathfrak p, V/\mathfrak p)$$ 
where $\underline{\Der}$ denotes the sheaf of derivations. 
The map $\psi_p(V,\nabla)$ vanishes precisely when $(V/{\mathfrak p}V,\nabla)$ has a full set of solutions, 
i.e. when $V/{\mathfrak p}V$ is spanned by its subsheaf of horizontal sections.  

\begin{conjecture}[Grothendieck $p$-Curvature Conjecture]
\label{conjecture:grothendieck}
Suppose that for almost all maximal ideals $\mathfrak p$ of $R$, the $p$-curvature of $(V/{\mathfrak p}V,\nabla)$ vanishes.  
Then the (complex) differential equation $(V,\nabla)$ has a full set of algebraic solutions, i.e. it becomes trivial on a finite \'etale cover of $X$.
\end{conjecture}

The truth of the conjecture does not depend on the choice of $R$ (see \S IV of \cite{Ka3}).

Katz (\cite{Ka1}, Theorem 13.0) showed that the condition on $p$-curvatures implies that $(V,\nabla)$ has regular singular points.  
The Riemann-Hilbert correspondence then implies that the conclusion of Conjecture \ref{conjecture:grothendieck} is equivalent to asking 
that the monodromy representation $\rho:~\pi_1(X(\C))\to\GL(n,\C)$ associated to $(V,\nabla)$ have finite image.  
In the same paper, Katz showed that if the $p$-curvatures vanish, and $X \hookrightarrow \bar X$ is a smooth compactification 
with $\bar X \backslash X$ a normal crossings divisor, then the local monodromy of $(V,\nabla)$ around a divisor at infinity is finite. 


Katz also made the following more general conjecture, which describes the Lie algebra $\mathfrak g$ of the Tannakian Galois group $G_{\gal}$ of $(V,\nabla).$
Recall that if $K(X)$ denotes the function field of the $R$-scheme $X,$ then $G_{\gal}$ is an algebraic group over $K(X),$ 
and $\mathfrak g$ is a Lie algebra over $K(X).$ As above one can ``reduce'' $\mathfrak g$ modulo almost any prime of $R.$ We then have 

\begin{conjecture}[Katz]
\label{conjecture:katz}
With the notation of Conjecture \ref{conjecture:grothendieck}, $\mathfrak g$ is the smallest algebraic Lie subalgebra of 
$\End_{\mathcal O_X} V \otimes K(X)$ such that 
for almost every prime $\mathfrak p$ of $R,$ the reduction of $\mathfrak g$ modulo $\mathfrak p$ contains the image of  
$\psi_p(V/{\mathfrak p}V,\nabla)$.
\end{conjecture}

Let $\rho:\pi_1(X(\C)) \rightarrow \GL_n(\C)$ be the monodromy representation of $(V,\nabla).$
When $(V,\nabla)$ has regular singular points, $\mathfrak g$ is essentially (up to change of base ring) the Lie algebra 
of the Zariski closure of $\rho(\pi_1(X(\C))).$ Thus Katz's conjecture describes the topological monodromy in terms 
of $p$-curvatures.

Katz proved that for any $X,$ Conjecture \ref{conjecture:grothendieck} is true for all $(V,\nabla)$ on $X$ 
if and only if Conjecture \ref{conjecture:katz} is true for all $(V,\nabla)$ on $X.$  For any particular $(V,\nabla)$ one knows only the obvious implication 
Conjecture \ref{conjecture:katz} $\implies$ Conjecture \ref{conjecture:grothendieck}.

Yet another result of Katz \cite{Ka2} says that Conjecture \ref{conjecture:grothendieck} is true when $(V,\nabla)$ is a Picard-Fuchs equation 
\footnote{That is, $(V,\nabla)$ arises from the de Rham cohomology of a smooth $X$-scheme.} 
or a suitable direct factor of one. Under some conditions, this result was extended by Andr\'e to $(V,\nabla)$ 
in the Tannakian category generated by a Picard-Fuchs equation [An, \S 16]. 

Following work of the Chudnovkys, Andr\'e [An] and then Bost [Bo] proved both conjectures in the case when the monodromy representation is solvable. 
So far there have been no examples of varieties $X$ with nonabelian fundamental group, for which the full conjectures are known. 
In this paper we give such examples; namely, we prove that the conjectures hold for many locally symmetric varieties.  

Recall that the set of complex points of a locally symmetric variety has the form 
$M=\Gamma\backslash G(\R)^+/K$, where $G$ is a connected, 
reductive, adjoint algebraic group over $\Q$, 
$G(\R)^+ \subset G(\R)$ denotes the connected component of the identity, $K \subset G(\R)$ 
is a maximal compact subgroup of $G$ with $G(\R)/K$ a Hermitian symmetric domain, 
and $\Gamma \subset G(\Q) \cap G(\R)^+$ is an arithmetic subgroup. 

Note that since $G$ is an adjoint group, we have 
$G = G_1 \times \dots \times G_m$ with each $G_i$ adjoint and $\Q$-simple. 

\begin{theorem}
\label{theorem:locsym}
Let $M$  be a locally symmetric variety, as above, and suppose that each factor $G_i$ has $\Q$-rank $\geq 1,$ and $\R$-rank $\geq 2.$
Then the Grothendieck-Katz conjecture holds for any $(V,\nabla)$ over $M$.
\end{theorem}

Note that Theorem \ref{theorem:locsym} omits the case when $M$ is a curve, a case which would imply (by standard methods) 
the full Conjecture \ref{conjecture:katz}. The proof of the theorem combines Katz's result on local monodromy with the Margulis normal subgroup theorem.

Symmetric varieties are naturally equipped with a collection of variations of Hodge structure (see \S 3 below for a more precise discussion). 
The underlying vector bundles with connection of these variations of Hodge structure form an interesting class of examples to which 
Conjecture \ref{theorem:locsym} may be applied. In this case the Lie algebra $\mathfrak g$ in Conjecture \ref{conjecture:katz} is a factor of $\Lie G.$ 
We refer to these bundles as {\it automorphic bundles}.

Using Margulis superrigidity, one can show that essentially any $(V,\nabla)$ on $M$ with regular singular points is a direct 
summand of an automorphic bundle. Combining this with Andr\'e's extension of Katz's result on Picard-Fuchs equations, 
allows us to  give another proof of Conjecture \ref{conjecture:katz} in the case when $G$ has no $\Q$-simple factors 
of exceptional type. 

\section{Proof of the main theorem}

\subsection{Toroidal compactifications}

As above, let $G$ be a connected, reductive, adjoint group over $\Q,$ and let $K \subset G(\R)^+$ be a maximal compact subgroup such that 
$D=G(\R)^+/K$ is a Hermitian symmetric domain.  Let $\Gamma$ be a torsion-free arithmetic subgroup in $G(\R)^+ \cap G(\Q),$ and let 
$M=\Gamma\backslash G(\R)^+/K$ be the corresponding locally symmetric variety.  

We recall the theory of toroidal compactifications for locally symmetric varieties \cite{AMRT}. 
Such a compactification $\bar{M}^{tor} \supset M$ is a projective variety admitting $M$ as a dense open subset. 
The identity map extends to a (unique) continuous map from $\bar{M}^{tor}$ to the Baily-Borel compactification. 
By construction the boundary  $\bar{M}^{tor} \backslash M$
is a union of non-empty closed subsets indexed by maximal parabolic subgroups $G$ which are defined over $\Q.$ 

The toroidal compactification is not unique; but depends on a choice of combinatorial data.  
This data can always be chosen so that the complement $\bar{M}^{\tor}$ is a divisor with normal crossings. 
A loop in $M$ around an irreducible component of  $\bar{M}^{\tor} \backslash M$ corresponds 
to a non-trivial unipotent element in $\Gamma$ (see \cite{AMRT}, III, \S 5, Main Thm.~I).


\subsection{Proof of Theorem \ref{theorem:locsym}}

First we may replace $\Gamma$ by a finite index subgroup, and assume that $\Gamma$ has the form 
$\Gamma_1 \times \dots \times \Gamma_m$ where $\Gamma_i \subset G_i(\Q)\cap G_i(\R)^+$ is an arithmetic subgroup. 
To prove the theorem it suffices to prove the analogous result for each $G_i.$ Replacing $G$ by $G_i,$ 
we may assume that $G$ is $\Q$-simple. 

By the result of Katz (\cite{Ka3}, Thm.~10.2) already mentioned above, it suffices to prove Conjecture \ref{conjecture:grothendieck} 
for an arbitrary vector bundle with integrable connection $(V,\nabla)$ on $M.$ Thus we may assume that the $p$-curvatures 
of $(V,\nabla)$ vanish for almost all $\mathfrak p.$

Now consider a smooth toroidal compactification $\bar M^{\tor}$ of $M,$ such that the boundary 
$\bar M^{\tor}\backslash M$ is a normal crossing divisor. 
Since the $\Q$-rank of $G$ is $\geq 1,$ $G$ has a parabolic subgroup defined over $\Q$ and 
the boundary is nonempty. Let $\gamma \in \Gamma$ be an element represented by a loop around 
an irreducible component of $\bar M^{\tor}\backslash M.$ Then $\gamma$ is a nontrivial unipotent 
element, and in particular has infinite order.


According to \cite{Ka1}, Thm.~13.0, since the $p$-curvatures of $(V,\nabla)$ are almost all $0,$ the local monodromy groups of $(V,\nabla)$ around the components of $\bar{M}^{\tor}\backslash M$ 
are finite groups. In particular $\rho(\gamma)$ has finite order and the kernel of $\rho$ is infinite.  

Now  since the $\Q$-rank of $G$ is $\geq 1,$ $G(\R)$ has no compact factors. 
Since $G$ is $\Q$-simple, $\Gamma \subset G(\R)$ is an irreducible lattice in a real semisimple Lie group of real rank at least 
$2.$ These properties of $G$ and $\Gamma$ imply that we can apply the Margulis Normal Subgroup Theorem (\cite{Ma}, IX, 6.14), which states that any normal 
subgroup of such a lattice is finite and central or has finite index.  Since $\ker(\rho) \vartriangleleft \Gamma$ is infinite, 
it follows that $\ker(\rho)$ has finite index in $\Gamma$, i.e. that $\rho(\Gamma)$ is finite, and we are done.
\qed

\section{Automorphic bundles and superrigidity}

\subsection{Automorphic bundles}

Let $G,$ $M$ and $\Gamma$ be as in the introduction, and 
let $\tilde G$ be a connected reductive group over $\Q$ with $\tilde G^{\ad} = G.$ 
Recall that there is a unique $G(\R)$ conjugacy class of 
cocharacters $h: \C^\times \rightarrow G(\R)$ such that the centralizer 
of $h$ in $G(\R)$ is compact. Fix such a cocharacter, and suppose that $h$ 
lifts to $\tilde h: \C^\times \rightarrow \tilde G(\R).$

Attached to any representation of $\tilde G$ (defined over $\Q$) there is a canonical 
variation of $\Q$-Hodge structure on $M$ [De, 1.1.12]. 
In particular, a representation of $\tilde G$ gives rise to an algebraic vector 
bundle on $M$ equipped with an integrable connection, having regular singular points. 
We call a vector bundle with connection arising in this way an {\it automorphic bundle}. 

For example, if $M$ is (a component of) a moduli space of principally polarized abelian varieties, 
and $\tilde G$ is a general symplectic group equipped 
with its standard representation, then the corresponding automorphic bundle arises from degree $1$ cohomology of 
the universal family of abelian varieties over $M.$ 

\begin{theorem}
\label{theorem:classification} Let $G$ be $\Q$-simple with $\R$-rank $\geq 2,$ and $\Q$-rank $\geq 1,$ 
and $(V,\nabla)$ a vector bundle with connection on $M$ having regular singular points. 
Then there exists a finite covering $\tilde M \rightarrow M$ such that the restriction of 
$(V,\nabla)$ to $\tilde M$ is a direct summand of an automorphic bundle.
\end{theorem}
\proof We will write $G_{\R}$ for $G$ viewed as an algebraic group over $\R.$ 

Since $(V,\nabla)$ has regular singular points it suffices to show that its associated 
monodromy representation is a direct summand of one attached to an automorphic bundle.
Let $\rho: \Gamma \rightarrow \GL_n(\C)$ denote the monodromy representation of $(V,\nabla).$ 
Viewing $\C$ as an $\R$-vector space, and replacing $n$ by $2n,$ 
we may assume that $\rho$ factors through $\GL_n(\R).$

Now let $\tilde G^{\der}$ denote the universal cover of $G.$ After replacing $\Gamma$  
by a finite index subgroup, we may view $\Gamma$ as a subgroup of $\tilde G^{\der}.$ 
Moreover, we may assume that $\rho$ factors through $\SL_n(\R).$ 
By Margulis superrigidity [Wi, \S12.A, Ex 8], after again replacing $\Gamma$ by a finite index subgroup, 
$\rho_{\alg}$ is induced by a map of algebraic $\R$-groups $\rho_{\alg}: \tilde G^{\der}_{\R} \rightarrow \SL_n.$

Let $Z = \ker(\tilde G^{\der} \rightarrow G)$ be the center of $\tilde G^{\der}.$ 
Choose a $\Q$-torus $T,$ and an embedding of $\Q$-groups $Z \hookrightarrow T.$ 
Let $K/\Q$ be a quadratic imaginary extension, and set 
$$ \tilde G = \tilde G^{\der} \times_Z R_{K/\Q} T.$$ 
where $R_{K/\Q}$ denotes restriction of scalars from $K$ to $\Q.$  
Then $h: R_{\C/\R} \G_m \rightarrow G_{\R}$ lifts to $\tilde G.$ 
Fix such a lifting. 

Let $\tilde G \hookrightarrow \GL(W_{\Q})$ 
be a faithful representation on a $\Q$-vector space $W_{\Q}.$ 
Then $W_{\R} = W_{\Q}\otimes_{\Q}\R$ is a faithful 
representation of $\tilde G^{\der}_{\R},$ so $\rho_{\alg}$ is a direct summand of 
a $\tilde G^{\der}_{\R}$-representation of the 
form $W_{\R}^{\otimes n}\otimes W_{\R}^{*\otimes m},$ where $W_{\R}^*$ 
denotes the dual of $W_{\R}.$
Then the composite 
$$ \Gamma \hookrightarrow \tilde G^{\der}(\R) \overset {\rho_{\alg}} \rightarrow \SL_n(\R). $$ 
is the direct summand of $W_{\R}^{\otimes n}\otimes W_{\R}^{*\otimes m},$ whose complexification 
is a representation associated to an automorphic bundle. 
\qed

We can use the argument in the above proof together with Andr\'e extension of Katz's result 
on Picard-Fuchs equation to give another proof of the main theorem when $G$ has no factors 
of exceptional type. 

Andr\'e's theorem says that Conjecture \ref{conjecture:katz} holds when $(V,\nabla)$ is in the Tannakian category 
generated by a Picard-Fuchs equation attached to a smooth projective morphism $f:Y \rightarrow X,$ 
{\it provided} one of the fibres $Y_x$ of $f$ has a connected motivic Galois group. 
Andr\'e defines this group using his theory of motivated cycles. 
Here we remark only that when $f$ is a family of abelian varieties, the motivic Galois group is equal to 
the Mumford-Tate group of $Y_x,$ and hence is connected [An, 16.2, 16.3].   

\begin{theorem} Suppose that $G$ is $\Q$-simple with $\R$-rank $\geq 2,$ and $\Q$-rank $\geq 1,$ and assume 
that $G$ is of type $A,B,C,$ or $D.$ Then the Grothendieck-Katz conjecture 
holds for any $(V,\nabla)$ on $M.$
\end{theorem}

\proof Let $\tilde G^{\der}$ be as in the proof of the previous theorem. Since we may replace $\Gamma$ 
by a finite index subgroup we may assume that $\Gamma \subset \tilde G^{\der}(\R),$ as before. 

Let $(V,\nabla)$ on $M$ have vanishing $p$-curvatures for almost all $\mathfrak p,$ and $\rho$ the corresponding 
representation of $\Gamma.$ As above, we may assume that $\rho$ is induced by a representation $\rho_{\alg}$ of 
the real algebraic group $\tilde G^{\der}.$ Let $V_{\C}$ denote the underlying complex 
vector space of $\rho_{\alg},$ and decompose $V_{\C}$ as a sum of representations 
on each of which the center $Z \subset \tilde G^{\der}$ acts via a character. It suffices consider each summand 
individually, so we may assume that $\rho(Z(\C))$ consists of scalars.

For any $g \geq 1,$ let $\mathcal H_g$ denote $\Sp_{2g}(\R)$ modulo a maximal compact subgroup. 
The condition on the type of the $\Q$-simple factors of $G$ implies (\cite{De} \S 1.3) 
that for some $g \geq 1,$ there exists an arithmetic subgroup $\Gamma_g \subset \Sp_{2g}(\R),$ 
a finite index subgroup $\Gamma' \subset \Gamma,$ and a map of complex varieties 
$$ \Gamma'\backslash G(\R)/K \rightarrow \Gamma_g\backslash {\mathcal H}_g $$ 
such that the corresponding map on fundamental groups $\Gamma' \rightarrow \Gamma_g $ is induced 
by a map of $\Q$-groups $\tilde G^{\der} \rightarrow \Sp_{2g}$ with finite central kernel. 

Let $W_{\Q}$ denote the representation of $\tilde G^{\der}$ which is induced by the standard representation 
of $\Sp_{2g}.$ Then $W_{\Q}$ is attached to an automorphic bundle $(W,\nabla_W)$ which arises as the cohomology 
of a family of abelian varieties on $M.$ Hence by Andr\'e's result, mentioned above, 
the Grothendieck-Katz conjecture holds for all bundles in the Tannakian category generated by $(W,\nabla_W).$ 

Now the bundle $(\End V, \nabla_{\End V})$ has regular singular points, and the corresponding monodromy 
representation $\End \rho$ factors through $G,$ as $Z$ acts on $V_{\C}$ via scalars. Hence $\End V$ is in the 
Tannakian category generated by $(W, \nabla_W).$ It follows that $(\End \rho)(\Gamma)$ is a finite group, 
and replacing $\Gamma$ by a finite index subgroup we may assume it is trivial.

It follows that $\Gamma$ acts on $V_{\C}$ by scalars, and hence so does $\tilde G^{\der},$ as $\Gamma$ is Zariski 
dense in $\tilde G^{\der}.$ As the latter group is connected and semisimple, we conclude that $\Gamma$ acts trivially 
on $V_{\C}.$
\qed

\noindent
Dept. of Mathematics, University of Chicago\\
5734 University Ave.\\
Chicago, Il 60637\\


\begin{thebibliography}{AMRT}
\small

\bibitem[AMRT]{AMRT}
A.~Ash, D.~Mumford, M.~Rapoport, Y.~Tai, 
Smooth compactifications of locally symmetric varieties, 
{\it Lie Groups: History, frontiers and applications IV} (1975). 

\bibitem[An]{An}
Y.~Andr\'e, Sur la conjecture des $p$-courbures de Grothendieck-Katz et un probl\`eme de Dwork, 
{\it Geometric aspects of Dwork theory}, 55-112, (2004). 

\bibitem[Bo]{Bo}
J.B. Bost, Algebraic leaves of algebraic foliations over number fields, {\it Publ. Math. }IHES, Vol. 93 (2001), 161-221.

\bibitem[De]{De}
P.~Delgine, Vari\'et\'es de Shimura: interpr\'etation modulaire et techniques de construction de mod\`eles canoniques, 
{\it Automorphic forms, representations and $L$-functions (Corvallis 1977)}, Proc.~Sym.~Pure Math. 33 (1979), 247-290. 

\bibitem[Ka1]{Ka1}
N.~Katz,  Nilpotent connections and the monodromy theorem: Applications of a result of Turrittin, 
{\it Publ. Math} IHES,  No. 39 (1970), 175--232. 

\bibitem[Ka2]{Ka2}
N.~Katz,  Algebraic solutions of differential equations ($p$-curvature and the Hodge filtration), {\it Invent. Math.} 18 (1972), 1--118.

\bibitem[Ka3]{Ka3} 
N.~Katz,  A conjecture in the arithmetic theory of differential equations, {\it Bull.~SMF} 110 (1982), 203-239. 

\bibitem[Ma]{Ma}
G.~Margulis, Discrete subgroups of semisimple Lie groups, 
{\it Ergebnisse der Math.} 17 (1991).

\bibitem[Wi]{Wi}
D.~Witte, Introduction to Arithmetic Groups,  
{\it preliminary version}, 2008.

\end{thebibliography}
\end{document}